\documentclass[11pt,a4paper]{article}

\usepackage{a4wide}
\usepackage{graphicx}
\usepackage{latexsym}
\usepackage{epsfig}
\usepackage{amssymb}
\usepackage{amstext}
\usepackage{amsgen}
\usepackage{amsxtra}
\usepackage{amsgen}
\usepackage{amsthm}

\newtheorem{thm}{Theorem}[section]
\newtheorem{prop}[thm]{Proposition}

\newtheorem{cor}[thm]{Corollary}
\newtheorem{rem}[thm]{Remark}

\theoremstyle{definition}

\newtheorem{assumption}[thm]{Assumption}

\numberwithin{equation}{section}

\newcommand{\dive}{\nabla\cdot}
\newcommand{\e}{\varepsilon}
\newcommand{\ee}{^{\e}}
\newcommand{\zz}{^{0}}

\newcommand{\R}{\mathbb{R}}
\newcommand{\T}{\mathbb{T}}

\newcommand{\intt}{\int_{0}^{t}}

\newcommand{\rhot}{\widetilde{\rho}}
\newcommand{\vt}{\widetilde{v}}
\newcommand{\ct}{\widetilde{c}}
\newcommand{\rhob}{\bar \rho}
\newcommand{\vb}{\bar v}
\newcommand{\cb}{\bar c}
\newcommand{\rhou}{\underline{\rho}}
\newcommand{\vu}{\underline{v}}
\newcommand{\cu}{\underline{c}}
\newcommand{\rhoc}{\widehat{\rho}}
\newcommand{\vc}{\widehat{v}}

\newcommand{\Ut}{\widetilde{U}}
\newcommand{\Ub}{\bar U}
\newcommand{\Uc}{\widehat{U}}
\newcommand{\Uu}{\underline{U}}
\newcommand{\E}{\mathcal{E}}
\newcommand{\intto}{\int_{\T^2}}

\begin{document}

\title{Singular convergence of nonlinear hyperbolic chemotaxis systems to Keller--Segel type models}

\author{M. Di Francesco and D. Donatelli}

\maketitle

\begin{abstract}
In this paper we deal with diffusive relaxation limits of
nonlinear systems of Euler type modeling chemotactic movement of
cells toward Keller--Segel type systems. The approximating systems
are either hyperbolic--parabolic or hyperbolic--elliptic. They all
feature a nonlinear pressure term arising from a \emph{volume
filling effect} which takes into account the fact that cells do
not interpenetrate. The main convergence result relies on
compensated compactness tools and is obtained for large initial
data under suitable assumptions on the approximating solutions. In
order to justify such assumptions, we also prove an existence
result for initial data which are small perturbation of a constant
state. Such result is proven via classical Friedrichs's
symmetrization and linearization. In order to simplify the
coverage, we restrict to the two--dimensional case with periodical
boundary conditions.
\end{abstract}

\section{Introduction}

This paper deals with diffusive relaxation limits for the
nonlinear hyperbolic model describing chemotactic movement of
cells, also known as the \emph{persistence and chemotaxis} model,
\begin{equation}\label{eq:beforescaling_intro}
  \begin{cases}
        \partial_{\tau}\rho +\dive (\rho v)=0 & \\
        \partial_{\tau} v + v\cdot \nabla v + \nabla
        g(\rho)=\nabla c - d v & \\
        \sigma\partial_{\tau} c = \Delta c + \alpha \rho -\beta c,
  \end{cases}
\end{equation}
with $\alpha, \beta, d,\sigma$ positive constants. The function
$g(\rho)$ is taken of the form $g(\rho)=\rho^\gamma$ with
$\gamma>0$, we shall discuss this choice later on in this
introduction. The model (\ref{eq:beforescaling_intro}) has been
introduced and motivated very precisely in
\cite{ambrosi_gamba_serini}, whereas similar models have been also
discussed in
\cite{preziosi_et_al,kowalczyk_gamba_preziosi,dolak_hillen,filbet_laurencot_perthame}.
We shall briefly summarize the biological motivations behind
(\ref{eq:beforescaling_intro}) by framing them in the general
context of PDE systems describing chemotactical phenomena.

The analysis of partial differential equations modeling chemotaxis
goes back to the work of Keller and Segel \cite{KS70}, who
proposed a macroscopic model for aggregation of cellular slime
molds, and to the earlier related work of Patlak \cite{Pat53}, who
derived similar models with applications to the study of
long-chain polymers. In the successive decades, the term
\emph{chemotaxis} has been used to represent the dynamics of
several biological systems (such as the bacteria Escherichia Coli,
or the amoebae Dyctiostelium Discoideum, or endothelial cells of
the human body responding to angiogenic factors secreted by a
tumor) in which the motion of a species is biased by the gradient
of a certain chemical substance. Typically these models consist of
a system of drift--diffusion equations of the form
\begin{equation}\label{eq:limit_intro}
\begin{cases}
  \rho_t = \Delta \rho -\dive \left(\rho\chi(\rho,c)\nabla
  c\right) & \\
  c_t =\Delta c +r(\rho,c), &
  \end{cases}
\end{equation}
with diffusion terms modeling random motion for the density $\rho$
of the individuals (cells, bacteria and so on) and for the
concentration $c$ of the \emph{chemoattractant} (the chemical
substance responsible of the chemotactical movement), first order
drift terms modeling chemotactical aggregation and zero order
reaction terms in the equation for the chemoattractant. The
coefficient $\chi(\rho,c)$ is called \emph{chemotactical
sensitivity}. The simplest model combining diffusion and
chemotaxis is the well known parabolic--elliptic
Patlak--Keller--Segel system (or simply Keller--Segel system)
\begin{equation}\label{KeSeintro}
\begin{cases}
  \rho_t = \Delta \rho -\dive \left(\chi \rho\nabla
  c\right) & \\
  0=\Delta c +\rho, &
  \end{cases}
\end{equation}
which has the interesting mathematical feature of producing smooth
solutions for small initial norms (in the appropriate space) and
blow-up (in the form of concentration to deltas) for large initial
norms. The rigorous analysis of such mathematical issue (also
extended to fully parabolic systems and to more complex models)
has attracted the attention of many mathematicians in the last
decades. We mention the pioneering works of J\"{a}ger--Luckaus
\cite{JL92}, Nagai \cite{Nag95}, Herrero--Velazquez \cite{HV96}
among others. The existence vs. blow--up problem in two space
dimensions for the classical Keller--Segel model (\ref{KeSeintro})
has been completely solved in \cite{DP04}, where the authors
proved that if the initial mass is less than a threshold value
$m^*$ (depending on the coefficient $\chi$) then the solution
exists globally in time in $L^1$, whereas if the initial mass is
larger than $m^*$ then the solution blows up in a finite time. We
refer to the surveys \cite{HorstmannI,HorstmannII} for a complete
and detailed description of the literature of this topic.

In the last years, some authors
\cite{kowalczyk,carrillo_calvez,LS06} have proposed variants of
Keller--Segel type models featuring global existence of solutions
no matter how large the initial mass is, obtained by replacing the
linear diffusion term in the equation for the population density
by a \emph{degenerate nonlinear diffusion} term with super--linear
growth for large densities. This choice can be motivated by taking
into account the fact that cells do not interpenetrate (that is,
they are full bodies with nonzero volume) and therefore diffusion
is supposed to inhibit singular aggregation effects when the
density is very high. We mention here that other authors proposed
the use of a nonlinear chemotaxis coefficient $\rho\chi(\rho)$
which attains the value zero when the population density $\rho$
reaches a fixed maximal value -- see for instance
\cite{PaHil,BO04,BDiFD06} -- being this choice motivated by the
fact that individuals stop aggregating when the density is too
high. In both cases, in the resultant model overcrowding of cells
(concentration to deltas for the cells density $\rho$) is
prevented independently on any initial parameter.

In the last years, several authors have started to describe
biological systems with chemotaxis from a hydrodynamical point of
view, i. e. via nonlinear hyperbolic systems of Euler type, see in
particular
\cite{ambrosi_gamba_serini,preziosi_et_al,kowalczyk_gamba_preziosi,dolak_hillen,filbet_laurencot_perthame}.
The models obtained are of the form of system
(\ref{eq:beforescaling_intro}), where the chemotactical force
$\nabla c$ in (\ref{eq:beforescaling_intro}) and the pressure
contribute to balance the rate of change of the momentum.
Moreover, our model (\ref{eq:beforescaling_intro}) features a
friction term modeling the drag between cells and the substrate
material (some authors also considered models with a linear
viscous term). In this framework, the nonlinear pressure term
$g(\rho)$ in (\ref{eq:beforescaling_intro}) plays the role of the
diffusive one in the drift--diffusion equation. Therefore, one can
interpret the overcrowding--preventing effect described before
(sometimes referred to as \emph{volume filling effect}) by
thinking of the cellular matter as a medium with \emph{limited
compressibility}, i. e. closely packed cells exhibit a limited
amount of resistance to compression. In this sense, a reasonable
choice of a pressure $g(\rho)$ is a function of the form
$g(\rho)=\rho^\gamma$, $\gamma>0$. Such an expression also has the
advantage of modeling absence of stresses for low densities (see
\cite{ambrosi_gamba_serini} for a more detailed description).

In this paper we want to contribute to the problem of establishing
a rigorous mathematical link between the system
(\ref{eq:beforescaling_intro}) and several Keller Segel type
models of the form (\ref{eq:limit_intro}) in terms of diffusive
relaxation limits. A typical example of diffusive scaling on the
system (\ref{eq:beforescaling_intro}) that we shall consider (see
subsection \ref{sec:firstscaling}) is the following
\begin{equation*}
  d =\frac{1}{\epsilon},\qquad \tau=\frac{t}{\e},\quad
  v^\e(x,t)=\frac{1}{\e}v\left(x,\frac{t}{\e}\right)
\end{equation*}
which transforms (\ref{eq:beforescaling_intro}) into the following
rescaled system
\begin{equation}\label{eq:afterscalingintro}
\begin{cases}
        \partial_{t}\rho\ee +\dive (\rho\ee v\ee)=0 & \\
        \e^2\left[\partial_{t} v\ee + v\ee\cdot \nabla v\ee\right] + \nabla
        g(\rho\ee)=\nabla c\ee - v\ee & \\
        \e\partial_{t} c\ee = \Delta c\ee + \alpha \rho\ee -\beta
        c\ee.
  \end{cases}
\end{equation}
Formally, as $\e\rightarrow 0$, we expect the solution
$(\rho^\e,v^\e,c^\e)$ to system (\ref{eq:afterscalingintro}) to
behave like the solution $(\rho\zz,u\zz,c\zz)$ to the
drift--diffusion system of Keller--Segel type
\begin{equation}\label{eq:reducedintro}
\begin{cases}
\partial_{t}\rho\zz+\dive (\rho\zz \nabla (c\zz-g(\rho\zz) )=0 & \\
\Delta c\zz +\alpha \rho\zz -\beta c\zz=0, &
\end{cases}
\end{equation}
where the loss of the persistence term in the equation for the
momentum yields a constitutive law for the velocity $v\zz=\nabla
c\zz- \nabla g(\rho\zz)$ (which can be considered as an equivalent
of the Darcy law in \cite{MM90}).

A way to understand the meaning of this phenomena is to consider
it as the large time behaviour of dissipative nonlinear hyperbolic
systems and to look at the asymptotic profile as the relaxed
equilibrium. This is the case for many relevant situations in
mathematical physics and applied mathematics. Singular limits with
a structure similar to \eqref{eq:afterscalingintro} have been
analyzed by Marcati and Milani \cite{MM90}. In that paper they
investigate the porous media flow as the limit of the Euler
equation in $1-D$, later generalized by Marcati and Rubino
\cite{MR00} to the $2-D$ case. Relaxation phenomena of the same
nature appear also in the zero relaxation limits for the
Euler-Poisson model for semiconductors devices and they were
investigated by Marcati and Natalini \cite{MN95a, MN95b} in the
$1-D$ case and by Lattanzio and Marcati \cite{lattanzio_marcati}
in the multi-D case. For a general overview of the theory of the
singular limits see the survey \cite{DM02} and the paper
\cite{DM04}, where the theory is completely set up.

To perform the relaxation limit we follow the same techniques
developed in \cite{MM90,MR00,DM04} (among others), which are
crucially based on the method of compensated compactness of Tartar
and Murat (see \cite{Tar79,Tar83,Mur78}) combined with the Young
measures associated to the relaxing sequence $\rho\ee$ (see
\cite{Tar79, DiP83a, DiP83b, DiP85b, DiP85a}). Throughout the
whole paper, we shall restrict ourselves to the case of two space
dimensions, which is also the most treated case in the literature
concerning Keller--Segel type systems. Moreover, for the sake of
simplicity we shall work on the $2$-dimensional torus.
 We shall prove that this singular limit can be rigorously justified as far as the new time
variable $\tau$ stays in a bounded interval $[0,T]$ for an
arbitrary $T$ and provided that certain a priori assumptions holds
for the solution to (\ref{eq:afterscalingintro}) (see assumption
\ref{basicassumption} below). These a priori assumptions are usual
in the framework of relaxation limits for nonlinear hyperbolic
systems (see also \cite{lattanzio_marcati}, \cite{MR00}) and they don't include
any smallness assumption on the initial conditions. The rigorous statements of these results are contained in Theorem \ref{tconv}.

In order to produce a nontrivial class of solutions to the
nonlinear hyperbolic system (\ref{eq:beforescaling_intro}) which
relax toward a Keller--Segel type model after a proper rescaling,
we shall also provide an existence theorem for the approximating
system (\ref{eq:afterscalingintro}) and prove the uniform
estimates needed to justify the assumptions
(\ref{basicassumption}) in case of initial densities $\rho_0$
which are small perturbation of an arbitrary non zero constant
state (see Theorem \ref{tpert}). This result is achieved by means
of the classical FRiedrichs' symmetrization technique and by a
linearization argument, see \cite{Fri54,KM81}. We remark that, in
many of the estimates performed here, the pressure term need not
to be of the form $g(\rho)=\rho^\gamma$. Indeed, some of the
estimates proven are still valid if one considers a logarithmic
pressure $g(\rho)=\log \rho$, which corresponds to a linear
diffusion term in the limit problem (\ref{eq:limit_intro}) (this
fact is not in contradiction with the blow--up of the density in
the limit problem with linear diffusion, see the Remark
\ref{remblowup}). However, while considering the alternative
scaling introduced in subsection \ref{sec:scalingpoisson} (where
the limit is the classical Keller--Segel system
(\ref{KeSeintro})), such an expression for the pressure seems to
be essential in order to achieve the needed estimates no matter
how large the initial mass is, in a very similar fashion to what
happens in \cite{carrillo_calvez}. We remark that our convergence
results hold on an arbitrary time interval. Therefore, at least in
the case of the second scaling treated in \ref{sec:scalingpoisson}
(where the expression $g(\rho)=\rho^\gamma$ is crucial in order to
achieve the desired estimates), our result can be seen as a new
interpretation of the overcrowding--preventing effect due to the
power--like expression of the pressure. More precisely, the global
smoothness of the limit density $\rho^0$ (and the absence of
concentration to deltas for all times of $\rho^0$ as a byproduct)
can be obtained as a consequence of our relaxation result,
alternatively to the more direct proof developed in
\cite{kowalczyk,carrillo_calvez}.

The paper is organized as follows. In chapter \ref{chap:prelim} we
state the three different scalings we shall deal with. In chapter
\ref{chapest} we perform the uniform estimate we need in order to
prove the main convergence theorem. In chapter \ref{chapconv} we
prove the main convergence theorem for large data under the a
priori assumption \ref{basicassumption} by means of compensated
compactness and Minty's argument. In chapter \ref{chappert} we
prove an existence theorem for the approximating rescaled system
in order to provide a class of solutions satisfying the basic
assumptions (\ref{basicassumption}).

\section{Preliminaries and rescalings}\label{chap:prelim}

We consider the following \emph{persistence and chemotaxis} model
\begin{equation}\label{eq:beforescaling}
  \begin{cases}
        \partial_{\tau}\rho +\dive (\rho v)=0 & \\
        \partial_{\tau} v + v\cdot \nabla v + \nabla
        g(\rho)=\nabla c - d v & \\
        \sigma\partial_{\tau} c = \Delta c + \alpha \rho -\beta c
  \end{cases}
\end{equation}
where $\alpha, \beta, d,\sigma$ are nonnegative constants. The
system \eqref{eq:beforescaling} is endowed with the following
$1$--periodic initial data
\begin{equation*}
\rho(0,x)=\rho_{0}(x), \qquad v(0,x)=v_{0}(x), \qquad c(0,x)=c_{0}(x)
\end{equation*}
The nonlinear function $g(\rho)$ has the form
$$g(\rho)=\rho^{\gamma},\ \hbox{for some}\
\gamma>0.$$
\begin{rem}
\emph{The nonlinear function $g(\rho)$ grows faster than
$\kappa^{\ast}\log \rho$, for large $\rho$, where
$\kappa^{\ast}=M/4\pi$ and $M$ is the total mass of $\rho$. More
precisely, there exists $\mathcal{U}>0$, such that
$$\text{for any $\rho\geq \mathcal{U}$ and $\kappa>\kappa^{\ast}$}\quad g(\rho)\geq \kappa\log\rho,$$
(see \cite{carrillo_calvez}). \label{r1}}
\end{rem}
Some of the results contained in the present paper hold in any
space dimension $n$, whereas some of them are true only in the
case $n=2$. In order to simplify the coverage, we shall always
restrict ourselves to the latter case. In the sections
\ref{chap:prelim}, \ref{chapest} and \ref{chapconv} we shall not
deal with the existence theory of (\ref{eq:beforescaling}),
whereas we shall work under the following basic assumption.
\begin{assumption}\label{basicassumption}
There exists a global solution $(\rho,v,c)$ to
(\ref{eq:beforescaling}), smooth enough in order to justify the
estimates contained in section \ref{chapest} and such that
\begin{itemize}
  \item [(A1)] the total mass $M=\int \rho dx$ is conserved,
  \item [(A2)] $\rho(x,t)\geq k>0$
  \item [(A3)] $(\rho, \rho v)\in L^\infty(\T^2 \times [0,+\infty))$.
\end{itemize}
\end{assumption}

Let us now explain in detail the relaxation limits we want to
perform. We shall deal with three different asymptotic regimes for
(\ref{eq:beforescaling}), corresponding to small parameter limits
of three different types of scaling.

\subsection{First scaling: large times and large
damping}\label{sec:firstscaling}

For a fixed constant $\e>0$ we consider the large damping rate
$d =\frac{1}{\epsilon}$ in (\ref{eq:beforescaling}). Then, we
introduce the fast time variable $$\tau=\frac{t}{\e},$$ and the
new independent variables
\begin{equation}
v^{\e}(x,t)=\frac{1}{\e}v(x,\tau),\quad
\rho^{\e}(x,t)=\rho(x,\tau),\quad c^{\e}(x,t)=c(x,\tau).
\label{s1}
\end{equation}
Moreover, we fix $\sigma=1$ in the third equation. Then, system
(\ref{eq:beforescaling}) in the new variables reads
\begin{equation}\label{eq:afterscaling}
\begin{cases}
        \partial_{t}\rho\ee +\dive (\rho\ee v\ee)=0 & \\
        \e^2\left[\partial_{t} v\ee + v\ee\cdot \nabla v\ee\right] + \nabla
        g(\rho\ee)=\nabla c\ee - v\ee & \\
        \e\partial_{t} c\ee = \Delta c\ee + \alpha \rho\ee -\beta
        c\ee.
  \end{cases}
\end{equation}
The formal limit as $\e\rightarrow 0$ is given by the
parabolic--elliptic system
\begin{equation}\label{eq:reduced}
\begin{cases}
\rho\zz_t+\dive (\rho\zz \nabla (c\zz-g(\rho\zz) )=0 & \\
\Delta c\zz +\alpha \rho\zz -\beta c\zz=0. &
\end{cases}
\end{equation}

\subsection{Second scaling: large time and large damping with Poisson
coupling}\label{sec:scalingpoisson}

A simplified version of (\ref{eq:beforescaling}), namely with
$\beta=0$ and $\sigma=0$, is given by the following system
\begin{equation}\label{eq:beforescaling2}
  \begin{cases}
        \partial_{\tau}\rho +\dive (\rho v)=0 & \\
        \partial_{\tau} v + v\cdot \nabla v + \nabla
        g(\rho)=\nabla c - \gamma v & \\
        0= \Delta c + \alpha \rho.
  \end{cases}
\end{equation}
By performing the same scaling as before, namely
\begin{equation}
\tau=\frac{t}{\e},\quad v^{\e}(x,t)=\frac{1}{\e}v(x,\tau),\quad
\rho^{\e}(x,t)=\rho(x,\tau),\quad c^{\e}(x,t)=c(x,\tau),
\label{s2}
\end{equation}
 and by
putting $d=\frac{1}{\epsilon}$, we obtain
\begin{equation}\label{eq:afterscaling2}
\begin{cases}
        \partial_{t}\rho\ee +\dive (\rho\ee v\ee)=0 & \\
        \e^2\left[\partial_{t} v\ee + v\ee\cdot \nabla v\ee\right] + \nabla
        g(\rho\ee)=\nabla c\ee - v\ee & \\
       0= \Delta c\ee + \alpha \rho\ee.
  \end{cases}
\end{equation}

The formal limit as $\e\rightarrow 0$ leads to the usual
Keller--Segel model with nonlinear diffusion (see
\cite{carrillo_calvez})
\begin{equation}\label{eq:reduced2}
\begin{cases}
\rho\zz_t+\dive (\rho\zz \nabla (c\zz-g(\rho\zz) )=0 & \\
\Delta c\zz +\alpha \rho\zz =0. &
\end{cases}
\end{equation}

\subsection{Third scaling: diffusive scaling with small reaction
rates}\label{sec:thirdscaling}

Starting once again by (\ref{eq:beforescaling}), we consider the
case $\sigma=d=1$. We consider $\epsilon$--depending reaction
coefficients $\alpha$ and $\beta$, namely we require
$\alpha=\epsilon\widetilde{\alpha}$ and
$\beta=\epsilon\widetilde{\beta}$ for fixed $\widetilde{\alpha},
\widetilde{\beta}>0$. We then perform the diffusive scaling
\begin{equation}
\tau=\frac{t}{\epsilon^2},\quad y=\frac{x}{\epsilon},\quad v^{\e}(x,t)=\frac{1}{\e}v(x,\tau),\quad
\rho^{\e}(x,t)=\rho(x,\tau),\quad c^{\e}(x,t)=c(x,\tau).
\label{s3}
\end{equation}
This
leads to the rescaled system
\begin{equation}\label{eq:afterscaling3}
\begin{cases}
        \partial_{t}\rho\ee +\dive (\rho\ee v\ee)=0 & \\
        \e^2\left[\partial_{t} v\ee + v\ee\cdot \nabla v\ee\right] + \nabla
        g(\rho\ee)=\nabla c\ee - v\ee & \\
       \partial_t c\ee= \Delta c\ee + \widetilde{\alpha} \rho\ee
       -\widetilde{\beta}c.
  \end{cases}
\end{equation}
Therefore, the formal limit as $\epsilon\rightarrow 0$ is given in
this case by the following fully parabolic model (we drop the \
$\widetilde{}$\, symbol for simplicity)
\begin{equation}\label{eq:reduced3}
\begin{cases}
\rho\zz_t+\dive (\rho\zz \nabla (c\zz-g(\rho\zz) )=0 & \\
c\zz _t =\Delta c\zz +\alpha \rho\zz -\beta c\zz.&
\end{cases}
\end{equation}
\begin{rem}
\label{r2} \emph{From the hypotheses (A3) and the scalings
\eqref{s1}, \eqref{s2}, \eqref{s3}, it follows that the sequences
$\{\rho\ee\}$, $\{\e\rho\ee v\ee\}$ are uniformly bounded in
$L^{\infty}(\T^{2}\times [0,+\infty))$ with respect to $\e$.}
\end{rem}

\section{Estimates}\label{chapest}

In this section we provide suitable estimates on the solutions of
the three rescaled models (\ref{eq:afterscaling}),
(\ref{eq:afterscaling2}) and (\ref{eq:afterscaling3}). For future
use we define
\begin{equation}\label{generalizedpressure}
  P(\rho):=\int_{0}^{\rho}g(n)dn=\frac{1}{\gamma+1}\rho^{\gamma+1}.
\end{equation}

\subsection{First scaling}

We have the following (standard) energy estimate for the rescaled
system (\ref{eq:afterscaling}).
\begin{prop}\label{p1} The following identity is
satisfied for any $t\in[0,T]$, by any solution $(\rho\ee,v\ee,c\ee)$ to
(\ref{eq:afterscaling}):
\begin{align}
  & \int_{\T^{2}} \left[\frac{\e^2}{2}\rho\ee(x,t)|v\ee(x,t)|^2
  +P(\rho\ee(x,t))\right]dx+\frac{1}{2}\intt\int_{\T^{2}}\rho\ee(x,s)|v\ee(x,s)|^2dxds\nonumber
  \\
  & \  =\int_{\T^{2}} \left[\frac{\e^2}{2}\rho_{0}\ee(x)|v_{0}\ee(x)|^2
  +P(\rho_{0}\ee(x))\right]dx+\left(\tilde{K}t+\e\int_{\T^{2}}\frac{c_{0}\ee(x)}{2}dx\right).
  \label{estimate1}
\end{align}
\end{prop}

\proof By multiplying  second equation in
(\ref{eq:afterscaling}) by $\rho\ee v\ee$ by using the first equation  in (\ref{eq:afterscaling}) and by integration by parts it follows that
\begin{align}
&\frac{d}{dt}\int_{\T^{2}}\left[\frac{\e^2}{2}\rho\ee(x,t)|v\ee(x,t)|^2
  +P(\rho\ee(x,t))\right]dx+\int_{\T^{2}}\rho\ee(x,s)|v\ee(x,s)|^2 dx=\notag\\
&\int_{\T^{2}}\rho\ee(x,s)v\ee(x,s)\nabla c\ee(x,t)dx\leq \frac{1}{2}\int_{\T^{2}}\rho\ee(x,s)|v\ee(x,s)|^2 dx+\frac{1}{2}\|\rho\ee\|_{\infty}\int_{\T^{2}}|\nabla c\ee(x,t)|^{2}dx.
 \label{estimate2}
  \end{align}
Now, by multiplying the third equation of \eqref{eq:afterscaling} by $c\ee$ we get for any $\delta>0$
\begin{align}
\frac{d}{dt}\int_{\T^{2}}\frac{\e}{2}|c\ee(x,t)|^{2}dx&=-\int_{\T^{2}}|\nabla c\ee(x,t)|^{2}dx+\alpha\|\rho\ee\|_{\infty}\left(\frac{|\T^{2}|^{2}}{4\delta}+\delta\int_{\T^{2}}|c\ee(x,t)|^{2}dx\right)\notag\\
&-\beta\int_{\T^{2}}|c\ee(x,t)|^{2}dx
 \label{estimate3}
\end{align}
By choosing $\delta<\frac{\beta}{2}$, by integrating in time we obtain, for fixed constant $\tilde{K}$, independent on $\e$, that $c\ee$ satisfies the following inequality
\begin{equation}
\int_{\T^{2}}\frac{\e}{2}|c\ee(x,t)|^{2}dx+\frac{\beta}{2}\int_{0}^{t}\int_{\T^{2}}|c\ee(x,t)|^{2}dxds+\int_{0}^{t}
\int_{\T^{2}}|\nabla c\ee(x,t)|^{2}dxds\leq \tilde{K} t+\e\int_{\T^{2}}\frac{|c_{0}\ee(x)|^{2}}{2}dx.
\label{estimate4}
\end{equation}
The estimate \eqref{estimate1} follows now by using together \eqref{estimate2} with \eqref{estimate4} and by taking into account the hypothesis (A3).
\endproof
\begin{cor}
\label{c1}
Let $(\rho\ee,v\ee,c\ee)$ be a solution to (\ref{eq:afterscaling})
satisfying assumption \ref{basicassumption} and with initial datum
$(\rho\ee_0,v\ee_0,c\ee_0)$ satisfying
\begin{equation}\label{i1}
  \int_{\T^{2}}\left[\frac{\e^2}{2}\rho\ee_0|v\ee_0|^2dx + (\rho\ee_0)^{\gamma+1}+\e|c\ee_0|^2\right]dx\quad\hbox{uniformly bounded w.r.t.}\
  \e\ll1.
\end{equation}
Then, for all $T>0$,
\begin{align}
  &\e\sqrt{\rho\ee}v\ee &\quad&\text{is uniformly bounded in
  $L^{\infty}([0,T],L^p(\T^{2}))$, for all $p\geq 1$,}\label{c1.1}\\
&\rho\ee &\quad& \text{is uniformly bounded in
  $L^{\infty}([0,T],L^{p}(\T^{2}))$, for all $p\geq 1$,}\label{c1.2}\\
  &\sqrt{\rho\ee}v\ee &\quad& \text{is uniformly bounded in
  $L^{2}([0,T]\times\T^{2})$,}\label{c1.3}\\
%&\nabla c\ee &\quad&\text{ is uniformly bounded in
  %$L^{\infty}([0,T]\times \T^{2})$.}\label{c1.4}\\
&\sqrt{\e} c\ee &\quad&\text{ is uniformly bounded in
  $L^{\infty}([0,T],L^{2}(\T^{2}))$,}\label{c1.5}\\
& c\ee &\quad&\text{ is uniformly bounded in
  $L^{2}([0,T],H^{1}(\T^{2}))$.}\label{c1.6}
\end{align}
\end{cor}
\proof \eqref{c1.1} and \eqref{c1.2} are a consequence of the
assumption \ref{basicassumption}, while \eqref{c1.3} follows from
the inequality \eqref{estimate1}. Finally
 \eqref{c1.5}, \eqref{c1.6} follow from \eqref{estimate4}.
\endproof

\subsection{Second scaling}

We consider the following energy for the solution to
(\ref{eq:afterscaling2})
\begin{equation}\label{energypoisson}
  E_\e(t)=\int_{\R^2}\left[\frac{\e^{2}}{2}\rho\ee |v\ee|^2 + P(\rho\ee)
  -\frac{1}{2}\rho\ee c\ee\right]dx,
\end{equation}
where $P$ is given by (\ref{generalizedpressure}). For semplicity we will take $\alpha=1$. We have the following estimate.

\begin{prop}
\label{p2}
The following inequality is valid for a solution
$(\rho\ee,v\ee,c\ee)$ to (\ref{eq:afterscaling2}):
\begin{equation}
 E_{\e}(t)+\int_0^t\int_{\T^{2}}\rho\ee(x,s)|v\ee(x,s)|^2 dxds= E_{\e}(0)
 \label{e1}
 \end{equation}
\end{prop}

\proof By using the Poisson equation of the system \eqref{eq:afterscaling2} we
easily have
\begin{align*}
    \frac{d}{dt} E_\e(t) &= -\int_{\T^{2}} \rho\ee |v\ee|^2 dx
    +\int_{\T^{2}}\rho\ee v\ee\cdot\nabla c\ee dx +\int_{\T^{2}}c\ee
    \Delta c_{t}\ee dx\\
    &=-\int_{\T^{2}} \rho\ee |v\ee|^2 dx+\int_{\T^{2}}\rho\ee v\ee\cdot\nabla c\ee dx+\int_{\T^{2}}(c\ee\nabla\cdot(\rho\ee v\ee))dx,
    \end{align*}
and this implies
$$ \frac{d}{dt} E_\e(t)=-\int_{\T^{2}}\rho\ee|v\ee|^{2}dx.$$
Integration with respect to time yields \eqref{e1}.
\endproof
In order to recover an estimate for $\nabla c\ee$, let us introduce the following convex functional
$$J[\rho\ee]=\int_{\T^{2}}(P(\rho\ee(t))-\rho\ee c\ee(t))dx.$$
Now we proceed by estimating the functional $J[\rho\ee]$ from
below, using the same strategy of \cite{carrillo_calvez}. Let us
recall that if $c\ee\in W^{1,1}(\T^{2})$, then the convex
functional $J[\rho\ee]$ has a critical point $\rho^{\ast}$ which
is a solution of
\begin{equation}
g(\rho^{\ast})-c\ee=\lambda
\label{l1}
\end{equation}
whenever $\rho^{\ast}>0$ and  null otherwise. Here $\lambda$ is the Lagrange multiplier associated to the constraint given by the mass conservation $\int \rho^{\ast}=M$ and fixed by this condition. We refer to \cite{carrillo_calvez} and (\cite{CJM01}, Proposition 5) for details. Therefore, we have
$$J[\rho\ee]\geq\int_{\T^{2}}(P(\rho^{\ast})-\rho^{\ast}c\ee)dx=\int_{\{\rho^{\ast}>0\}}(P(\rho^{\ast})-\rho^{\ast}g(\rho^{\ast})+\lambda \rho^{\ast})dx.$$
By taking into account the Remark \ref{r1} we can introduce the corrective term $R$ such that $g(\rho^{\ast})=\kappa\log\rho^{\ast}+R(\rho^{\ast})$, then we have
\begin{equation}
\label{l2}
J[\rho\ee]\geq\int_{\T^{2}}(P(\rho^{\ast})-\kappa\rho^{\ast}\log\rho^{\ast})dx-\int_{\{\rho^{\ast}>0\}}\rho^{\ast}R(\rho^{\ast})dx+\lambda M.
\end{equation}
Now, \eqref{l1} implies $\kappa\log \rho^{\ast}+R(\rho^{\ast})=\lambda+c\ee$, whenever $\rho^{\ast}>0$ and thus
$$\int_{\{\rho^{\ast}>0\}}exp\left(\frac{R(\rho^{\ast})}{\kappa}\right)\rho^{\ast}dx=e^{\lambda/\kappa}\int_{\{\rho^{\ast}>0\}}exp\left(\frac{c\ee}{\kappa}\right)dx,$$
so we have
\begin{equation}
\label{l3}
\lambda=\kappa\log\left(\int_{\{\rho^{\ast}>0\}}e^{R/\kappa}\rho^{\ast}dx\right)-\kappa\log\left(\int_{\{\rho^{\ast}>0\}}e^{c\ee/\kappa}dx\right).
\end{equation}
If we replace $\lambda$ by its expression in the inequality \eqref{l2}, we conclude that
\begin{align}
J[\rho\ee]&\geq\int_{\T^{2}}(P(\rho^{\ast})-\kappa\rho^{\ast}\log\rho^{\ast})dx-\int_{\{\rho^{\ast}>0\}}\rho^{\ast}R(\rho^{\ast})dx\notag\\
&+\kappa M\log\left(\int_{\{\rho^{\ast}>0\}}e^{R/\kappa}\rho^{\ast}dx\right)-\kappa M\log\left(\int_{\{\rho^{\ast}>0\}}e^{c\ee/\kappa}dx\right).
\label{l4}
\end{align}
By taking into account the Remark \ref{r1} we have that
$$\int_{\{\rho^{\ast}\geq \mathcal{U}\}}(P(\rho^{\ast})-\kappa\rho^{\ast}\log\rho^{\ast})dx\geq C.$$
On the other hand, we have
$$\int_{\{\rho^{\ast}< \mathcal{U}\}}(P(\rho^{\ast})-\kappa\rho^{\ast}\log\rho^{\ast})dx\geq -\left(\sup_{[0,\mathcal{U})}(P-\kappa\rho\log\rho)^{-}\right)|\T^{2}|.$$
Therefore,
$$\int_{\T^{2}}(P(\rho^{\ast})-\rho^{\ast}\log\rho^{\ast})dx$$
is uniformly bounded form below. Now, the Jensen inequality for the probability density $\rho^{\ast}/M$ over the set where $\rho^{\ast}>0$, gives us that
$$exp\int_{\{\rho^{\ast}>0\}}\left(\frac{R(\rho^{\ast})}{\kappa}\frac{\rho^{\ast}}{M}dx \right)\leq \int_{\{\rho^{\ast}>0\}} e^{R/\kappa}\frac{\rho^{\ast}}{M}dx,$$
ans thus
$$\kappa M\log\left(\int_{\{\rho^{\ast}>0\}} e^{R/\kappa}\frac{\rho^{\ast}}{M}dx\right)-
\int_{\{\rho^{\ast}>0\}}\rho^{\ast}R(\rho^{\ast})dx\geq 0. $$
Finally, we recall and use the Trudinger - Moser inequality \cite{Mos71, CY88, GZ98}.
\begin{thm}
Assume that $\Omega\subset \R^{2}$ is a $C^{2}$, bounded, connected domain. It exists a constant $C_{\Omega}$, such that for all $h\in H^{1}$ with $\int_{\Omega}h=0$ we have
$$\int_{\Omega}exp(|h|)dx\leq C_{\Omega}exp\left(\frac{1}{8\pi}\int_{\Omega}|\nabla h |^{2}dx\right).$$
\end{thm}
By applying the previous theorem to $c\ee/\kappa$ we obtain
$$\int_{\{\rho^{\ast}>0\}}e^{c\ee/\kappa}dx\leq \int_{\T^{2}}e^{c\ee/\kappa}dx\leq exp\left(\frac{1}{8\pi\kappa^{2}}\int_{\T^{2}}|\nabla c\ee|^{2}dx\right)$$
and thus
$$-\kappa M\log\left(\int_{\{\rho^{\ast}>0\}}e^{c\ee/\kappa}dx\right)\geq -\frac{M}{8\pi\kappa^{2}}\int_{\T^{2}}|\nabla c\ee|^{2}dx.$$
So by \eqref{l4} we have that
\begin{equation}
\label{l5}
J[\rho\ee]\geq C-\frac{M}{8\pi\kappa^{2}}\int_{\T^{2}}|\nabla c\ee|^{2}dx
\end{equation}
\begin{prop}
\label{p3}
Assume $(\rho\ee,v\ee,c\ee)$ be a solution to (\ref{eq:afterscaling2})
satisfying assumption \ref{basicassumption} then
\begin{equation}\label{l6}
\int_{\T^{2}}|\nabla c\ee|^{2}dx \qquad \text{is uniformly bounded}.
\end{equation}
\end{prop}
\proof
We can rewrite \eqref{e1} as
\begin{equation}
E_{\e}(0)=\int_{\T^{2}}\frac{\e^{2}}{2}\rho\ee |v\ee|^2dx+\int_{\T^{2}}J[\rho\ee]+\frac{1}{2}\int_{\T^{2}}|\nabla c\ee(t)|^{2} dx+\int_0^t\int_{\T^{2}}\rho\ee(x,s)|v\ee(x,s)|^2 dxds.
\label{l7}
\end{equation}
Combining \eqref{l7} with \eqref{l5} we get that
\begin{align}
E_{\e}(0)&\geq\int_{\T^{2}}\frac{\e^{2}}{2}\rho\ee |v\ee|^2dx+\int_0^t\int_{\T^{2}}\rho\ee(x,s)|v\ee(x,s)|^2 dxds\notag\\
&+C|\T^{2}|+\frac{1}{2}\left(1-\frac{M}{4\pi\kappa}\right)\int_{\T^{2}}|\nabla c\ee(t)|^{2} dx.
\label{l8}
\end{align}
Finally, Remark \ref{r1} implies $ \kappa>\kappa^{\ast}$, i.e. $\left(1-\frac{M}{4\pi\kappa}\right)>0$ and thus
$$\int_{\T^{2}}|\nabla c\ee|^{2}dx$$
 is uniformly bounded.
\endproof
\begin{cor}
\label{c2}
Let $(\rho\ee,v\ee,c\ee)$ be a solution to (\ref{eq:afterscaling2})
satisfying assumption \ref{basicassumption} and with initial datum
$(\rho\ee_0,v\ee_0,c\ee_0)$ satisfying
\begin{equation}\label{i2}
  \int_{\T^{2}}\left[\frac{\e^2}{2}\rho\ee_0|v\ee_0|^2dx + (\rho\ee_0)^{\gamma+1}-\frac{1}{2}\rho\ee_0 c\ee_0\right]dx\quad\hbox{uniformly bounded w.r.t.}\
  \e\ll1.
\end{equation}
Then, for all $T>0$,
\begin{align}
  &\e\sqrt{\rho\ee}v\ee &\quad&\text{is uniformly bounded in
  $L^{\infty}([0,T],L^p(\T^{2}))$, for all $p\geq 1$,}\label{c2.1}\\
&\rho\ee &\quad& \text{is uniformly bounded in
  $L^{\infty}([0,T],L^{p}(\T^{2}))$, for all $p\geq 1$,}\label{c2.2}\\
  &\sqrt{\rho\ee}v\ee &\quad& \text{is uniformly bounded in
  $L^{2}([0,T]\times\T^{2})$,}\label{c2.3}\\
& c\ee &\quad&\text{ is uniformly bounded in
  $L^{\infty}([0,T], H^{1}(\T^{2}))$.}\label{c2.4}
  \end{align}
\end{cor}
\proof \eqref{c2.4} follows from Proposition \ref{p3} and by
taking into account that we are in a periodic domain. \eqref{c2.1}
and \eqref{c2.2} are a consequence of the assumption
\ref{basicassumption}, while \eqref{c2.3} is a consequence of
\eqref{e1} and (\ref{l6}).
\endproof
\subsection{Third scaling}
With the same procedure as in the Proposition \ref{p1} we are able
to prove that
\begin{prop}
\label{p33}
Let $(\rho\ee,v\ee,c\ee)$ be a solution to (\ref{eq:afterscaling3})
satisfying assumption \ref{basicassumption} and with initial datum
$(\rho\ee_0,v\ee_0,c\ee_0)$ satisfying
\begin{equation}\label{i3}
  \int_{\T^{2}}\left[\frac{\e^2}{2}\rho\ee_0|v\ee_0|^2dx + (\rho\ee_0)^{\gamma+1}+|c\ee_0|^2\right]dx\quad\hbox{uniformly bounded w.r.t.}\
  \e\ll1.
\end{equation}
Then, for all $T>0$,
\begin{align}
  &\e\sqrt{\rho\ee}v\ee &\quad&\text{is uniformly bounded in
  $L^{\infty}([0,T],L^p(\T^{2}))$, for all $p\geq 1$,}\label{c3.1}\\
&\rho\ee &\quad& \text{is uniformly bounded in
  $L^{\infty}([0,T],L^{p}(\T^{2}))$, for all $p\geq 1$,}\label{c3.2}\\
  &\sqrt{\rho\ee}v\ee &\quad& \text{is uniformly bounded in
  $L^{2}([0,T]\times\T^{2})$,}\label{c3.3}\\
%&\nabla c\ee &\quad&\text{ is uniformly bounded in
%  $L^{\infty}([0,T]\times \T^{2})$.}\label{c3.4}\\
&c\ee &\quad&\text{ is uniformly bounded in
  $L^{\infty}([0,T],L^{2}(\T^{2}))\cap L^{2}([0,T],H^{1}(\T^{2}))$.}\label{c3.5}
  \end{align}
\end{prop}

\section{Strong convergence}\label{chapconv}

This section is devoted to the study of the relaxation of the systems
\eqref{eq:afterscaling}, \eqref{eq:afterscaling2}, \eqref{eq:afterscaling3} towards their formal limit \eqref{eq:reduced},  \eqref{eq:reduced2},  \eqref{eq:reduced3}, respectively. As a consequence of the Corollary \ref{c1} and the Propositions \ref{p3}, \ref{p33} we have that, extracting if necessary a subsequence,
$$\nabla c\ee\rightharpoonup \nabla c^{0}\quad \text{as $\e\downarrow 0$ weakly in $L^{2}([0,T]\times \T^{2})$ }.$$
This convergence for $c\ee$ is enough to pass into the limit in
\eqref{eq:afterscaling}, \eqref{eq:afterscaling2},
\eqref{eq:afterscaling3} to get in the sense of distribution
\eqref{eq:reduced},  \eqref{eq:reduced2},  \eqref{eq:reduced3},
respectively, p rovided that $\rho\ee$ converges in a strong
topology. In fact by the Remark \ref{r2} we know that
$\rho\ee\rightarrow \rho^{0}$ in $L^{\infty}$ $\ast-$weakly, while
by \eqref{c1.2}, \eqref{c2.2}, \eqref{c3.2} we have
$\rho\ee\rightharpoonup \rho^{0}$ weakly in $L^{p}$, for any
$p>1$. These convergence are clearly too weak to pass into the
limit in the nonlinear terms of \eqref{eq:afterscaling},
\eqref{eq:afterscaling2}, \eqref{eq:afterscaling3}. So, in this
section we will investigate the strong convergence of the
approximating sequence $\rho\ee$. The analysis of this convergence
reduces to the analysis of the convergence of quadratic forms with
constant coefficients via the classical compensated compactness
technique due to Tartar \cite{Tar79, Tar83} and Murat \cite{Mur78}
(see Dacorogna \cite{Dac82}). As we will see later on, these
techniques will apply in the same way to the three scalings
\eqref{s1}, \eqref{s2}, \eqref{s3}, so we will discuss them
together.  Let us recall the following theorem
\begin{thm}\label{cc}
{\bf (Tartar's Compensated compactness)}\\
Let us consider
\begin{enumerate}
\item
a bounded open set $\Omega\subset\mathbb{R}^{n}$;
\item
a sequence $\{l^{\nu}\}_{\nu = 1}^{\infty}$,
$l^{\nu}:~\Omega\subset \mathbb{R}^{n}
\longrightarrow\mathbb{R}^{m} $;
\item
a symmetric matrix
$\Theta:~\mathbb{R}^{m}\longrightarrow\mathbb{R}^{m} $;
\item
constants $a_{jk}^{i}\in \mathbb{R}$, $i = 1,\ldots,q$,
$j = 1,\ldots,m$, $k = 1,\ldots,n$.
\end{enumerate}
Let us define
\begin{align*}
f(\alpha) & = \left\langle \Theta \alpha,\alpha \right\rangle,
\quad\text{for all $\alpha \in \mathbb{R}^{m}$;} \\
\Lambda & = \left\{ \lambda \in \mathbb{R}^{m}:~
\exists \eta \in \mathbb{R}^{n}\setminus
\{0\}, ~\sum\limits _{j,k}a_{jk}^{i} \lambda_{j}
\eta_{k} = 0, i = 1,\ldots,q\right\}.
\end{align*}
Assume that
\begin{itemize}
\item[\bf(a)]
there exists $\widetilde{l}\in L^{2}\left(\Omega\right)$
such that $l^{\nu}\rightharpoonup \widetilde{l}$ in
$L^{2}\left(\Omega\right)$ as $\nu\uparrow\infty$;
\item[\bf(b)]
$\mathcal{A}^{i}l^{\nu} = \sum\limits _{j,k}a^{i}_{jk}\frac{\partial
l_{j}^{\nu}}{\partial x_{k}}$, $i=1,\ldots q$ are relatively compact in
$H^{-1}_{loc}\left(\Omega\right)$;
\item[\bf(c)]
$f_{|\Lambda}\equiv 0$;
\item[\bf(d)] there exists
$\widetilde{f}\in \mathbb{R}$ such that $f(l)\rightharpoonup
\widetilde{f}$ in the sense of measures $\mathcal{M}(\Omega)$.
\end{itemize}
Then we have $\widetilde{f} = f(\widetilde{l})$.
\end{thm}
\subsection{Weak convergence of $\rho\ee P(\rho\ee)$}
First of all we start by studying the weak convergence
 of  $\rho\ee P(\rho\ee)$. Our goal will be to prove that
 $$\rho\ee P(\rho\ee)\rightharpoonup \rho^{0}P(\rho^{0}),$$
where  $\rho^{0}$ is the weak limit of $\rho\ee$.   To this end  we are going to apply the Theorem \ref{cc} in the same spirit of \cite{MR00}. In order to fit the into the hypotheses of the Theorem \ref{cc} we rewrite the first two equations of the systems
\eqref{eq:afterscaling},  \eqref{eq:afterscaling2},  \eqref{eq:afterscaling3}, as
\begin{equation}
\begin{cases}
\rho\ee_{t} +m\ee_{x}+n\ee_{y} = 0 \\
\displaystyle
{\e^{2}m\ee_{t}+\left(\e^{2}\frac{(m\ee)^{2}}{\rho\ee}+\gamma P(\rho\ee)\right)_{x}+\left(\e^{2}\frac{m\ee n\ee}{\rho\ee}\right)_{y} = \rho\ee c\ee_{x}-m\ee}\\
\displaystyle{
\e^{2}n\ee_{t}+\left(\e^{2}\frac{m\ee n\ee}{\rho^{\ee}}\right)_{x}+\left(\e^{2}\frac{(n\ee)^{2}}{\rho\ee}+\gamma P(\rho\ee)\right)_{y} = \rho\ee c\ee_{y}-n\ee.}
\end{cases}
\label{4.1.2}
\end{equation}
where
\begin{equation}
v\ee=(v_{1}\ee, v_{2}\ee)\qquad \rho\ee v\ee=(\rho\ee v_{1}, \rho\ee v_{2})=(m\ee, n\ee).
\label{4.1.1}
\end{equation}
It will be usefull  rewrite  \eqref{4.1.2} in the following way
\begin{align}
\rho\ee_{t} +m\ee_{x}+n\ee_{y} &= 0 \notag\\
\gamma P(\rho\ee)_{x}&=-\e^{2}m\ee_{t}-\e^{2}\left(\frac{(m\ee)^{2}}{\rho\ee}\right)_{x}-\e^{2}\left(\frac{m\ee n\ee}{\rho\ee}\right)_{y} + \rho\ee c\ee_{x}-m\ee\label{4.1.11}\\
\gamma P(\rho\ee)_{y}&=-\e^{2}n\ee_{t}-\e^{2}\left(\frac{m\ee n\ee}{\rho\ee}\right)_{x}-\e^{2}\left(\frac{(n\ee)^{2}}{\rho\ee}\right)_{y}+\rho\ee c\ee_{y}-n\ee.\notag
\end{align}
By using \eqref{c1.3}, \eqref{c2.3}, \eqref{c3.3}, \eqref{l6},  and the assumption (A3) we get that $\rho\ee v\ee, \ \rho\ee\nabla c\ee\in L^{2}([0,T]\times \T^{2})$. In fact
\begin{align}
&\|\rho\ee v\ee\|_{L^{2}([0,T]\times \T^{2})}\leq \|\sqrt{\rho\ee}\|_{\infty}\|\sqrt{\rho\ee}v\ee\|_{L^{2}([0,T]\times \T^{2})}\label{cc1}\\
&\|\rho\ee\nabla c\ee\|_{L^{2}([0,T]\times \T^{2})}\leq \|\rho\ee\|_{\infty}\|\nabla c\ee\|_{L^{2}([0,T]\times \T^{2})}\label{cc2}
\end{align}
Moreover, by taking into account the assumptions (A2) and (A3) and \eqref{c1.3}, \eqref{c2.3}, \eqref{c3.3} we have  that $\e^{2}\left(\frac{(m\ee)^{2}}{\rho\ee}\right)_{x}$,  is relatively compact in $H^{-1}([0,T]\times \T^{2})$. In fact let us consider $\omega$ relatively compact in $[0,T]\times \T^{2}$, then by taking into account (A2), (A3) and the Remark \ref{r2} we have,
\begin{equation}
\left\|\e^{2}\left(\frac{(m\ee)^{2}}{\rho\ee}\right)_{x}\right\|_{H^{-1}(\omega)}\leq
\sup _{\|\phi\|_{H^{1}_{0}(\omega)}=1}\left|\left\langle
\e^{2}\left(\frac{(m\ee)^{2}}{\rho\ee}\right)_{x} ,
\phi\right\rangle\right|\leq \e\|\rho\ee
v\ee\|_{\infty}\frac{1}{\sqrt{k}}\|\sqrt{\rho\ee}v\ee\|_{L^{2}(\omega)}
\label{cc2bis}
\end{equation}
In a similar way it can be proved that $\e^{2}\left(\frac{m\ee
n\ee}{\rho\ee}\right)_{y}$, $\e^{2}\left(\frac{m\ee
n\ee}{\rho\ee}\right)_{x}$,
$\e^{2}\left(\frac{(n\ee)^{2}}{\rho\ee}\right)_{y}$,
$\e^{2}(\rho\ee v\ee)_{t}$ are relatively compact in
$H^{-1}([0,T]\times \T^{2})$. Now, \eqref{cc1}--\eqref{cc2bis},
imply that
\begin{equation}
\begin{pmatrix}
\rho\ee_{t}+m\ee_{x}+n\ee_{y}   \\
P(\rho\ee)_{x}   \\
P(\rho\ee)_{y}
\end{pmatrix}
\qquad
\text{is  relatively compact in $\left(H^{-1}_{loc}\right)^{3}$.}   \label{4.1.3}
\end{equation}
In order to fit into the framework of the Theorem \ref{cc} we set $x_{1}=x$, $x_{2}=y$, $x_{3}=t$, $l\ee=(m\ee, n\ee, \rho\ee, P(\rho\ee))$, hence $m=4$. In our case the differential constraints are $q=3$. So we can define the matrices $\mathcal{A}^{1}, \mathcal{A}^{2}, \mathcal{A}^{3}\in \mathcal{M}_{4\times 3}$, where $\mathcal{A}^{i}=\{a^{i}_{jk}\}$, $i=1,2,3$, $j=1, \ldots,4$, $k=1,2,3$ as follows:
\begin{equation*}
   \mathcal {A}^{1}=\begin{pmatrix}
   1 & 0 & 0 \\
        0 & 1 & 0 \\
0 & 0 & 0\\
0 & 0 & 0
\end{pmatrix}       \quad
\mathcal{A}^{2}=\begin{pmatrix}
0 & 0 & 0 \\
        0 & 0 & 0 \\
        0 & 0 & 0 \\
                1 & 0 & 0
\end{pmatrix}
\quad
\mathcal{A}^{3}=\begin{pmatrix}
0 & 0 & 0 \\
        0 & 0 & 0 \\
        0 & 0 & 0 \\
                0 & 1 & 0
 \end{pmatrix}.
 \end{equation*}
The characteristic manifold $\Lambda$ is then given by
\begin{equation*}
\Lambda =\left\{\lambda \in \R^{4} \mid \exists \xi \in
\R^{3}\setminus \{0\},  B(\xi,\lambda)=0\right\}
\end{equation*}
where
\begin{equation*}
B(\xi,\lambda)=\begin{pmatrix}
\lambda_{1}\xi_{1}+\lambda_{2}\xi_{2}+\lambda_{3}\xi_{3}\\
\lambda_{4}\xi_{1} \cr      \lambda_{4}\xi_{2}.
\end{pmatrix}
\end{equation*}
Therefore
\begin{equation*}
\Lambda =\left\{\lambda\in \R ^{4} \mid
det\begin{pmatrix}\lambda_{1} &
\lambda_{2} & \lambda_{3}   \\
\lambda_{4} & 0 & 0 \\
0 & \lambda_{4} & 0
\end{pmatrix}    \right\}=\left\{\lambda \in \R^{4} \mid
\lambda_{3}\lambda_{4}=0\right\}.
\end{equation*}
If we define
$$M=\frac{1} {2}\begin{pmatrix}
0 & 0 & 0 & 0\\
0 & 0 & 0 & 0\\
0 & 0 & 0 & 1\\
0 & 0 & 1 & 0
\end{pmatrix}    \in \mathcal {M}_{4\times 4}, $$
 then
$f(\lambda)={\lambda}^{T}M \lambda    =\lambda_{3}\lambda_{4}$ and, of course
$f_{|\Lambda}\equiv 0$.
Now, by applying the Theorem \ref{cc} we have
 $f(l\ee)\rightharpoonup f(\tilde{l})$, and in our case this means
\begin{equation*}
\rho\ee P(\rho\ee)\rightharpoonup \rho^{0} P^{0},
\end{equation*}
where $P^{0} =w-\lim P (\rho\ee)$.

\subsection{Strong convergence of $\rho\ee$}
In the previous section we proved that
\begin{equation*}
\rho\ee P(\rho\ee)\rightharpoonup \rho^{0} P^{0}.
\end{equation*}
Here we will be able to prove that
$$\rho\ee\rightarrow \rho^{0}\quad \text{strongly in $L^{p}_{loc}$, $p<+\infty$}.$$
At this point we can follow the methods of \cite{MM90}, \cite{MR00}.
First of all let as use  Minty's argument (\cite{L-JL84}, \cite{MM90}) to prove that
$P^0=P(\rho^0)$. Since the function $P$  is monotone,  for any $w\in
L^\infty$ e $\varphi$ test function, $\varphi >0$, we have that
\begin{equation}
H(\e)\equiv \int\!\!\int (P(\rho\ee)-P(w))(\rho\ee-w)\varphi dxdt\geq 0;
\label{5.1}
\end{equation}
But for $\e\downarrow 0$, we have that
$$\int\!\!\int P(\rho\ee)\rho\ee\varphi dxdt\rightarrow \int\!\!\int P^0
\rho^0\varphi dxdt,$$
So from \eqref{5.1}  we get that for  $\e \downarrow 0$
$$H(\e)\rightarrow H\equiv \int\!\!\int(P^0-p(w))(\rho^{0}-w) \varphi dxdt \geq
0$$
If we choose $w=\rho^{0} +\lambda z$, with $\lambda \leq 0$ and  arbitrary $z \in L^\infty$, we have
\begin{align*}
G(\lambda,z) &\equiv  \int\!\!\int (P^0-P(\rho^{0} +\lambda z))z\varphi dxdt\\
&
=\frac{1}{\lambda}\int\!\!\int (P^0-P(\rho^{0} +\lambda z))\lambda z\varphi dxdt
\leq 0
\end{align*}
and for  $\lambda \uparrow 0$, $G(0,z)\leq 0$ for  any $z \in L^\infty$,
then
$$G(0,z)=\int\!\!\int (P^0-P(\rho^{0}))z\varphi dxdt =0,$$
and finally  $P^0=P(\rho^0)$.\\
Our next step now, is to prove the strong convergence for
$\rho\ee \rightarrow \rho^0$ in $L^{p}_{loc}$. To this end we characterize the weak convergence by means of Young's probability measures
(see \cite{Tar79},\cite{DiP83a}, \cite{DiP83b}, \cite{DiP85b},
\cite{DiP85a}). Let us recall that if  $\{u\ee\}$ is sequence converging to
 $U$ in $L^\infty$ $\ast$-weakly , we can associate to
 $\{u\ee\}$ a family  $\{\nu_{(x,t)}(\lambda)\}$ of probability measures such that for any continuos function
 $F(\cdot)$
$$\ast-\lim_{\e \rightarrow 0} F(u\ee)=\int
F(\lambda)\nu_{(x,t)}(d\lambda)\qquad  a.e$$
If $\nu_{(x,t)}=\delta _{U(x,t)}$ then $u\ee \rightarrow U$ strongly in
$L^{p}_{loc}$ for any $p\in(1, +\infty)$ (see \cite{Dac82}, Corollary 6.2).
Let    $\{\nu_{(x,t)}\}$ be the  family of Young's probability measures
 associated  to the   sequence  $\{\rho\ee\}$: since
 $\rho\ee \longrightarrow \rho^{0}$ in $L^{\infty}$ $\ast$-weakly, we can find a closed interval  $\left[a,b\right]$, $0\leq a\leq b$,
such that  $supp \ \nu_{(x,t)}\subseteq \left[a,b\right]$. Since
$P(r) = r^{\alpha}$, $\alpha> 1$, we have three possibilities:
\begin{enumerate}
\item $P\in C^{2}\left(\R\setminus \{0\}\right)$ e $P''(r)\uparrow
+\infty$ for $r\downarrow 0$, if    $1< \alpha <2$;
\item $P\in C^{2}\left(\R\right)$ e $P''(0) = 1$, if            $\alpha = 2$;
\item $P\in C^{2}\left(\R\right)$ and $P''(0) = 0$, if          $\alpha >2$.
\end{enumerate}
Let us assume that $1< \alpha \leq 2$.  Then we can write for any $\lambda, \lambda_{0}$
$$P(\lambda) - P(\lambda_{0}) = P'(\lambda_{0})
\left(\lambda -     \lambda_{0}\right) + \frac {1}{2}
P''\left(\lambda^{*}\right)         \left(\lambda - \lambda_{0}\right)^{2},$$
where $\lambda^{*}$ belongs to the segment between   $\lambda$ and  $\lambda_{0}$.
If we choose
$$\lambda_{0} = \int_{a}^{b} \lambda    \nu_{(x,t)}\left(d\lambda\right)=\rho_{0},$$
since  $P^0=P(\rho^0)$
$$P\left(\lambda_{0}\right) = \int_{a}^{b}
P\left(\lambda\right)\nu_{(x,t)}\left(d\lambda\right)$$
so that
$$\int_{a}^{b} \{ P(\lambda) - P(\lambda_{0}) \}            \nu_{(x,t)}\left(d\lambda\right) = 0.$$
On the other hand we  also have
\begin{align*}
\int_{a}^{b}P'(\lambda_{0})\left(\lambda - \lambda_{0}\right)
\nu_{(x,t)} \left(d\lambda\right)
= P'(\lambda_{0})\left\{\int_{a}^{b} \lambda    \nu_{(x,t)}\left(d\lambda\right)-\lambda_{0}\int_{a}^{b}\nu_{(x,t)}\left(d\lambda\right)\right\}=0,
\end{align*}
so we can conclude that
$$\int_{a}^{b} P''\left(\lambda^{*}\right) \left(\lambda - \lambda_{0}
\right)^{2} \nu_{(x,t)}\left(d\lambda\right) = 0.$$
Taking
$E=\displaystyle{\min_{ \lambda \in \left[a,b\right]}}      P''\left(\lambda\right) > 0$, we get
$$E\int_{a}^{b} \left(\lambda - \lambda_{0}\right)^{2} \nu_{(x,t)}\left(d\lambda\right) \leq 0,$$
namely
$$\int_{a}^{b} \left(\lambda - \lambda_{0}\right)^{2}           \nu_{(x,t)}\left(d\lambda\right)= 0,$$
and it follows $a = b$ and $\nu_{(x,t)}=\delta$, a point mass and so we finally get
$$\rho\ee\rightarrow \rho_{0}\qquad \text{strongly in $L^{p}_{loc}$}.$$
To conclude we remark that in the case $\alpha>2$ this result can be obtained in the same way by using the function $-P^{-1}$.
\begin{rem}
\emph{The strong convergence result for $\rho\ee$ obtained in this
section is still valid in the case of linear diffusion, namely if
we consider $g(\rho)=\log \rho$ and consequently
$P(\rho)=\rho\log\rho -\rho$. }\end{rem} By using the estimates
and the strong convergence of the sequence $\{\rho\ee\}$ obtained
in the  previous section we are able to prove the following main
theorem.
\begin{thm}
\label{tconv} Let $T>0$ be arbitrary and let $(\rho\ee, v\ee,
c\ee)$ be a family of solutions to the system
\eqref{eq:afterscaling} (\eqref{eq:afterscaling2} and
\eqref{eq:afterscaling3} respectively) with initial data
satisfying (\ref{i1}) ((\ref{i2}) and (\ref{i3}) respectively).
Assume that the assumption \ref{basicassumption} holds, then,
there exist $\rho^{0}\in L^{\infty}([0,T]\times \T^{2})$ and
$c^{0}\in L^{2}([0,T],H^{1}(\T^{2}))$, such that, extracting if
necessary a subsequence,
\begin{align*}
\rho\ee&\longrightarrow\rho^{0} \qquad \text{strongly in $L^{p}_{loc}([0,T]\times \T^{2})$ for any $p<\infty$}\\
\nabla c\ee&\rightharpoonup \nabla c^{0} \qquad \text{weakly in $L^{2}([0,T]\times \T^{2})$.}
\end{align*}
Moreover the couple $(\rho^{0}, c^{0})$, satisfies the system
\eqref{eq:reduced} (\eqref{eq:reduced2} and \eqref{eq:reduced3}
respectively) in the sense of distributions.
\end{thm}

\section{Perturbation of constant states in the approximating
system} \label{chappert}

In this section we deal with the rescaled system
\begin{equation}\label{eq:pert1}
\begin{cases}
        \partial_{t}\rho\ee +\dive (\rho\ee v\ee)=0 & \\
        \e^2\left[\partial_{t} v\ee + v\ee\cdot \nabla v\ee\right] + \nabla
        g(\rho\ee)=\nabla c\ee - v\ee & \\
        \e\partial_{t} c\ee = \Delta c\ee + \alpha \rho\ee -\beta
        c\ee.
  \end{cases}
\end{equation}
with $x\in \T^2$, $t\geq 0$, where $\T^2$ is the flat normalized
two--dimensional torus. The system (\ref{eq:pert1}) is
complemented with the $1$-periodical initial data
\begin{equation*}
  \rho\ee(x,0)=\rho\ee_0(x),\quad v\ee(x,0)=v\ee_0(x),\quad
  c\ee(x,0)=c\ee_0(x).
\end{equation*}
We shall consider small perturbations of the stationary state
\begin{equation}\label{constantstates}
  (\rho,v,c)=(\rhot,\vt,\ct),\quad \rhot>0, \quad \vt=0,\quad
  \ct=\frac{\alpha}{\beta}\rhot
\end{equation}
and prove the existence of solutions $(\rho\ee,v\ee,c\ee)$ such
that the density $\rho\ee$ stays away from zero, uniformly in
$\varepsilon$, on a small enough time interval $[0,T]$ with $T$
independent on $\varepsilon$ (see similar results in \cite{KM81} and \cite{DiFM}). In order to perform this task, we
shall use an iterative method, namely we define recursively the
sequence $(\rho^n,v^n,c^n)$ as follows:
$(\rho^0(x,t),v^0(x,t),c^0(x,t))=(\rho\ee_0(x),v\ee_0(x),c\ee_0(x))$
and, for $n\geq 1$, $(\rho^n,v^n,c^n)$ solves the linear system
\begin{equation}\label{eq:sequence}
\begin{cases}
        \displaystyle{\partial_{t}\rho^n +v^{n-1}\cdot \nabla \rho^n + \rho^{n-1}\dive u^n=0} & \\
        \displaystyle{\partial_{t} v^n + v^{n-1}\cdot \nabla v^n + \frac{g'(\rho^{n-1})}{\e^2}\nabla
        \rho^n=\frac{1}{\e^2}\nabla c^n - \frac{1}{\e^2}v^n}\vspace{1mm} & \\
        \displaystyle{\partial_{t} c^n = \frac{1}{\e}\Delta c^n + \frac{\alpha}{\e} \rho^n
        -\frac{\beta}{\e}c^n.}
  \end{cases}
\end{equation}
From now on we shall drop the dependency on $\e$ on the solutions
$(\rho,v,c)$ to simplify the notation. Moreover, we shall use the
following notation: the variables taken at the step $n-1$ will be
denoted e. g. by $\rho^{n-1}=\rhoc$; the variables taken at the
step $n$ will be denoted without any further symbol, e. g. $\rho^n
= \rho$; the deviation from the aforementioned constant stationary
states will be denoted e. g. by $\rhob=\rho^n - \rhot$ and
$\rhou=\rho^{n-1}-\rhot$.

The first two equations in system (\ref{eq:sequence}) can be
easily viewed as a hyperbolic system in vectorial form. More
precisely, let us define the $3$--dimensional variable $U$ as
\begin{equation*}
  U:=(\rho,v^1,v^2),
\end{equation*}
where $v=(v^1,v^2)$. Let us denote
\begin{equation*}
  A_1(\Uc):=\left[\begin{matrix} \vc^1 & \rhoc & 0 \\
  \frac{g'(\rhoc)}{\e^2} & \vc^1 & 0 \\ 0 & 0 & \vc^1
  \end{matrix}\right],\qquad  A_2(\Uc):=\left[\begin{matrix} \vc^2 & 0 & \rhoc \\
  0 & \vc^2 & 0 \\ \frac{g'(\rhoc)}{\e^2} & 0 & \vc^2
  \end{matrix}\right],\qquad B(U):=\frac{1}{\e^2}\left(\begin{matrix} 0 \\
  \partial_{x_1}c - v^1 \\ \partial_{x_2}c -
  v^2\end{matrix}\right).
\end{equation*}
Then, with all these notations, the system (\ref{eq:sequence}) can
be rephrased as
\begin{equation}\label{eq:pert2}
\begin{cases}
\partial_t U + + A_1(\Uc) \partial_{x_1}U + A_2(\Uc)
\partial_{x_2}U= B(U) & \\
\displaystyle{\partial_t c = \frac{1}{\e}\Delta c
+\frac{\alpha}{\e}\rho - \frac{\beta}{\e}c}. &
\end{cases}
\end{equation}
The first line in (\ref{eq:pert2}) corresponds to a linear
hyperbolic system which can be \emph{symmetrized} by means of the
matrix
\begin{equation}\label{symmetrizer}
  S(\Uc):=\mathrm{diag}\left(\frac{g'(\rhoc)}{\e^2}, \rhoc, \rhoc
  \right).
\end{equation}
The matrix $S(\Uc)$ is uniformly positive definite provided the
variable $\rhoc$ satisfies a condition of the form $0<c\leq
\rhoc\leq C$ (we recall that $g'(\rho)=\gamma \rho^{\gamma-1}$
exhibits a singularity at zero in case of $\gamma<1$). The two
matrices $S(\Uc)A_1(\Uc)$ and $S(\Uc)A_2(\Uc)$ can be easily
proven to be symmetric. We now rewrite system (\ref{eq:pert2}) in
terms of the deviations $\Ub$ and $\cb$:
\begin{equation}\label{eq:pert3}
\begin{cases}
\partial_t \Ub + A_1(\Ut+\Uu) \partial_{x_1}\Ub + A_2(\Ut+\Uu)
\partial_{x_2}\Ub= B(\Ut+\Ub)=B(\Ub) & \\
\displaystyle{\partial_t \cb = \frac{1}{\e}\Delta \cb
+\frac{\alpha}{\e}\rhob - \frac{\beta}{\e}\cb}. &
\end{cases}
\end{equation}
We introduce the energy functional
\begin{equation*}
  \E(U,c):=\frac{1}{2}\intto \left[ U^T S(\Uc) U + \lambda c^2\right]dx =
  \frac{1}{2}\intto \left[ \frac{g'(\rhoc)}{\e^2}\rho^2 + \rhoc
  |v|^2 + \lambda c^2\right] dx,
\end{equation*}
where $\lambda>0$ is a constant to be chosen later on. We have the
following
\begin{prop} Let $T>0$. There exist constants $\e_0,\delta\in (0,1)$, $K\in(0,\rhot /2)$ such that, if
\begin{align}
  & \|\rho_0^\e-\rhot\|_{H^4(\mathbb{T}^2)} +
  \e\|v_0^\e\|_{H^4(\mathbb{T}^2)}+\sqrt{\e}\|c_0^\e-\ct\|_{H^4(\mathbb{T}^2)}
  \leq \delta \quad \hbox{and}\nonumber\\
  & \quad \sup_{0\leq t\leq T}\left(\|\rhou(t)\|_{H^4(\mathbb{T}^2)}+\e\|\vu(t)\|_{H^4(\mathbb{T}^2)} +\sqrt{\e}\|\cu(t)-\ct\|_{H^4(\mathbb{T}^2)}
  \right) \leq
  K,\label{ipotesi_induttiva}
\end{align}
for all $\e\in(0,\e_0)$, then,
\begin{equation}\label{final1}
  \sup_{0\leq t\leq
T}\left(\|\rhob(t)\|_{H^4(\mathbb{T}^2)}+\e
\|\vb(t)\|_{H^4(\mathbb{T}^2)}+\sqrt{\e}\|\cb(t)-\ct\|_{H^4(\mathbb{T}^2)}\right)
\leq K
\end{equation}
for all $\e\in(0,\e_0)$.
\end{prop}

\proof During the proof of this proposition we shall often make
use of the Sobolev inequality $\|f\|_{L^\infty(\T^2)}\leq
C\|f\|_{H^2(\T^2)}$.

\textsc{Step 1}. Due to the symmetry of the two matrices $SA_1$
and $SA_2$, we can use integration by parts in the evolution of
$\E(\Ub,\cb)$ as follows:
\begin{align}
  \frac{d}{dt} \E(\Ub,\cb) & = \frac{1}{2}\intto\Ub^T\left[S(\Uc)A_1(\Uc)\right]_{x_1}\Ub dx +
  \frac{1}{2}\intto\Ub^T\left[S(\Uc)A_2(\Uc)\right]_{x_2}\Ub dx\nonumber\\
   & \ + \intto\Ub^T S(\Uc)B(\Ub)dx -\frac{\lambda}{\e}\intto |\nabla
   \cb|^2 dx + \frac{\lambda\alpha}{\e}\intto\rhob \cb dx
   -\frac{\lambda\beta}{\e}\intto \cb^2 dx.\label{energyest1}
\end{align}
Due to the assumption (\ref{ipotesi_induttiva}) we have
\begin{align}
  \frac{d}{dt} \E(\Ub,\cb) & \leq
  C(K) \left(\|\nabla \rhou\|_{L^\infty} + \|\nabla \vu\|_{L^\infty}\right)
  \frac{1}{2}\intto\left(\frac{\rhob^2}{\e^2} + \frac{|\vb|^2}{\e^2}\right)
  dx+\frac{\|\rhoc\|_{L^\infty}}{\e^2}\intto \vb\cdot \nabla \cb dx\nonumber\\
  & \
  -\frac{(\rhot - K)}{\e^2}\intto |\vb|^2 dx-\frac{\lambda}{\e}\intto |\nabla
   \cb|^2 dx + \frac{\lambda\alpha}{\e}\intto\rhob \cb dx
   -\frac{\lambda\beta}{\e}\intto \cb^2 dx\label{energyest2}
\end{align}
for a function $C(K)>0$ of the constant $K$ such that $C$ is
continuous on $K\in [0,\rhot /2]$. By choosing
\begin{equation*}
  \lambda=\frac{(\rhot +K)^2}{\e(\rhot - K)}
\end{equation*}
and $\e_0< 1$, we can use once again (\ref{ipotesi_induttiva}) and
find a constant $C_1>0$ such that
\begin{align*}
     \frac{d}{dt}\E(\Ub,\cb)& \leq K C(K)\frac{1}{2}\intto\left(\frac{\rhob^2}{\e^2} + \frac{|\vb|^2}{\e^2}\right)
  dx -\frac{(\rhot-K)}{2\e^2}\intto |\vb|^2 dx -\frac{(\rhot
  +K)^2}{2(\rhot -K)
  \e^2}\intto|\nabla \cb|^2 dx\\
  & \ -\frac{(\rhot +K)^2\beta}{2(\rhot-K)\e}\intto
  \cb^2 dx + \frac{(\rhot+K)^2 \alpha^2}{2(\rhot-K)\e \beta}\intto \rhob^2
  dx.
\end{align*}
We now choose the constant $K$ such that $K C(K)<\frac{1}{2}(\rhot
-K)$. By using the coercivity property
\begin{equation}
  \E(U,c)\geq c(K)\intto\left[\frac{\rho^2}{\e^2} +|v|^2 +
  \frac{c^2}{\e}\right]dx,\label{coercivity}
\end{equation}
which holds for a certain $c(K)>0$, due to Gronwall inequality we
easily obtain
\begin{equation*}
\E(\Ub(t),\cb(t)) +\frac{1}{\e^2}\int_0^t\intto
\left[|\vb(\tau)|^2 dx + |\nabla \cb(\tau)|^2\right]dx d\tau +
\frac{1}{\e}\intt\intto \cb(\tau)^2 dxd\tau \leq
A\E(\Ub(0),\cb(0))e^{Bt}
\end{equation*}
for certain constants $A,B>0$ depending only on $K$ and $\e_0$.
The above implies in particular
\begin{equation*}
\sup_{0\leq t\leq T}\intto\left[\rhob(t)^2 +\e^2 |\vb(t)|^2 +\e
\cb(t)^2\right] dxdt \leq C(K) \delta e^{BT},
\end{equation*}
for a certain $C(K)$ depending on $K$. Therefore, by choosing
$\delta$ small enough such that $C(K)\delta e^{BT}\leq K^2$ we
obtain
\begin{equation*}
\sup_{0\leq t\leq
T}\left(\|\rhob(t)\|_{L^2(\mathbb{T}^2)}+\e\|\vb(t)\|_{L^2(\mathbb{T}^2)}
+\sqrt{\e}\|\cb(t)\|_{L^2(\mathbb{T}^2)}
  \right) \leq
  K.
\end{equation*}

\textsc{Step 2}. We now perform the energy estimate of the space
derivatives of $(\Ub,\cb)$. For $j=1,2$ we denote the derivatives
with respect to $x_j$ by the subscript
$\rho_j=\partial_{x_j}\rho$. The system satisfied by the first
derivatives of $(\Ub,\cb)$ is
\begin{equation}\label{eq:pert4}
\begin{cases}
\partial_t \Ub_j + A_1(\Ut+\Uu) \partial_{x_1}\Ub_j + A_2(\Ut+\Uu)
\partial_{x_2}\Ub_j= B(\Ub)_j- A_1(\Ut+\Uu)_j \Ub_1 - A_2(\Ut+\Uu)_j \Ub_2  & \\
\displaystyle{\partial_t \cb_j = \frac{1}{\e}\Delta \cb_j
+\frac{\alpha}{\e}\rhob_j - \frac{\beta}{\e}\cb_j}. &
\end{cases}
\end{equation}
The evaluation of the energy
\begin{equation*}
  \E(\Ub_j,\cb_j)=\frac{1}{2}\intto \left[ \Ub_j^T S(\Uc) \Ub_j + \lambda \cb_j^2\right]dx
\end{equation*}
in a similar way as in (\ref{energyest1}) yields
\begin{align*}
  \frac{d}{dt} \E(\Ub_j,\cb_j) & = \frac{1}{2}\intto\Ub_j^T\left[S(\Uc)A_1(\Uc)\right]_{x_1}\Ub_j dx +
  \frac{1}{2}\intto\Ub_j^T\left[S(\Uc)A_2(\Uc)\right]_{x_2}\Ub_j dx\\
   & \ + \intto\Ub_j^T S(\Uc)B(\Ub)_j dx -\frac{\lambda}{\e}\intto |\nabla
   \cb_j|^2 dx + \frac{\lambda\alpha}{\e}\intto\rhob_j \cb_j dx
   -\frac{\lambda\beta}{\e}\intto \cb_j^2 dx\\
   & \ -\intto \Ub_j^T S(\Uc)A_1(\Uc)_j \Ub_1 dx -\intto \Ub_j^T S(\Uc)A_2(\Uc)_j \Ub_2 dx.
\end{align*}
Assumption (\ref{ipotesi_induttiva}) allows for the estimate of
the first two terms above as in (\ref{energyest1}), as well as for
the estimate of the last two terms in a similar fashion. The
result is the following estimate
\begin{align*}
  \frac{d}{dt} \E(\Ub_j,\cb_j) & \leq
  \widetilde{C}(K) \left(\|\nabla \rhou\|_{L^\infty} + \|\nabla \vu\|_{L^\infty}\right)
  \frac{1}{2}\intto\left(\frac{\rhob_j^2}{\e^2} + \frac{|\vb_j|^2}{\e^2}\right)
  dx+\frac{\|\rhoc\|_{L^\infty}}{\e^2}\intto \vb_j\cdot \nabla \cb_j dx\\
  & \
  -\frac{(\rhot - K)}{\e^2}\intto |\vb_j|^2 dx-\frac{\lambda}{\e}\intto |\nabla
   \cb_j|^2 dx + \frac{\lambda\alpha}{\e}\intto\rhob_j \cb_j dx
   -\frac{\lambda\beta}{\e}\intto \cb_j^2 dx,
\end{align*}
which is the equivalent of the estimate (\ref{energyest2}) where
$(\Ub,\cb)$ are replaced by their first derivatives. Therefore, we
can easily conclude as before
\begin{equation*}
\sup_{0\leq t\leq
T}\left(\|\nabla\rhob(t)\|_{L^2(\mathbb{T}^2)}+\e\|\nabla\vb(t)\|_{L^2(\mathbb{T}^2)}
+\sqrt{\e}\|\nabla \cb(t)\|_{L^2(\mathbb{T}^2)}
  \right) \leq
  K.
\end{equation*}

\textsc{Step 3}. The second space derivatives of $(\Ub,\cb)$
satisfy the system
\begin{equation}\label{eq:pert5}
\begin{cases}
\partial_t \Ub_{ij} + A_1(\Ut+\Uu) \partial_{x_1}\Ub_{ij} + A_2(\Ut+\Uu)
\partial_{x_2}\Ub_{ij}&\hspace{-3mm}= B(\Ub)_{ij} -A_1(\Ut+\Uu)_i\Ub_{1j} - A_2(\Ut+\Uu)_i\Ub_{2j} \\
& \ - A_1(\Ut+\Uu)_{ij} \Ub_1 - A_2(\Ut+\Uu)_{ij} \Ub_2 \\
& \ -A_1(\Ut+\Uu)_j \Ub_{1j}-A_2(\Ut+\Uu)_j \Ub_{2j} \\
\displaystyle{\partial_t \cb_{ij} = \frac{1}{\e}\Delta \cb_{ij}
+\frac{\alpha}{\e}\rhob_{ij} - \frac{\beta}{\e}\cb_{ij}}, &
\end{cases}
\end{equation}
for $i,j=1,2$. The structure of system (\ref{eq:pert5}) is similar
to (\ref{eq:pert4}) and therefore the estimate of the energy
\begin{equation*}
  \E(\Ub_{ij},\cb_{ij})=\frac{1}{2}\intto \left[ \Ub_{ij}^T S(\Uc) \Ub_{ij} + \lambda \cb_{ij}^2\right]dx
\end{equation*}
can be performed as in step 2. The only extra terms which needs to
be analyzed are the following, for $i,j,k=1,2$ ($C(K)$ denotes a
generic constant depending on $K$):
\begin{align*}
   & \intto \Ub_{ij}^T S(\Uc)A_k(\Uc)_{ij}\Ub_k \leq \frac{C(K)}{\e^2} \intto |\Ub_{ij}|\left[|\Uu_i\Uu_j| +|\Uu_{ij}|\right]|\Ub_k|dx \\
   & \ \leq \frac{K^2 C(K)}{\e^2}\intto |\Ub_{ij}||\Ub_k|dx + \frac{K C(K)}{\e^2}\intto
   |\Ub_{ij}||\Uu_{ij}|dx\\
   & \ \leq \frac{(K+K^2) C(K)}{\e^2}\left[\intto |\Ub_{ij}|^2 dx
   + K^2\right],
\end{align*}
where we have used once again (\ref{ipotesi_induttiva}) and the
result in step 2. Notice that so far we have used $L^\infty$
estimates only up to the firs order derivatives of $\Ub$ and
$\Uu$. In the last inequality above, the second derivatives are
only estimated in $L^2$. We have therefore obtained, for $0<K<1$,
\begin{align*}
  \frac{d}{dt} \E(\Ub_{ij},\cb_{ij}) & \leq
  \frac{K^3 C(K)}{\e^2}+ C(K)K
  \frac{1}{2}\intto\left(\frac{\rhob_{ij}^2}{\e^2} + \frac{|\vb_{ij}|^2}{\e^2}\right)
  dx+\frac{\|\rhoc\|_{L^\infty}}{\e^2}\intto \vb_{ij}\cdot \nabla \cb_{ij} dx\\
  & \
  -\frac{(\rhot - K)}{\e^2}\intto |\vb_{ij}|^2 dx-\frac{\lambda}{\e}\intto |\nabla
   \cb_{ij}|^2 dx + \frac{\lambda\alpha}{\e}\intto\rhob_{ij} \cb_{ij} dx
   -\frac{\lambda\beta}{\e}\intto \cb_{ij}^2 dx
\end{align*}
and, by using the same choice of $\lambda$ and $K$ as in step 1,
after using Gronwall Lemma we obtain
\begin{equation*}
\E(\Ub(t),\cb(t)) \leq C(K)\left[\E(\Ub(0),\cb(0))
+\frac{K^3}{\e^2}\right] e^{Bt}.
\end{equation*}
Then, the coercivity property (\ref{coercivity}) and the
assumptions (\ref{ipotesi_induttiva}) imply
\begin{equation*}
\sup_{0\leq t\leq
T}\left(\|D^2\rhob(t)\|_{L^2(\mathbb{T}^2)}+\e\|D^2\vb(t)\|_{L^2(\mathbb{T}^2)}
+\sqrt{\e}\|D^2 \cb(t)\|_{L^2(\mathbb{T}^2)}
  \right) \leq C(K)(\delta + K^3)
\end{equation*}
and clearly, a choice of $\delta$ and $K$ small enough implies
$C(K)(\delta + K^3)<K^2$, which concludes the estimate of the
second derivatives.

\textsc{Step 4}. In order to conclude the proof of the
proposition, one needs to perform the same energy estimate also on
the space derivatives of order $3$ and $4$. All the estimates on
the nonlinear terms on the right--hand side are analogous to those
in Step 3. The integrals with over-quadratic terms always contains
not more than two terms involving more than two derivatives.
Therefore, all the extra terms can be estimated in $L^\infty$ by
using assumption (\ref{ipotesi_induttiva}) and the results in the
previous steps. We shall skip the details of these computations.
The proof is complete.
\endproof

We are now ready to state the main theorem of this section.

\begin{thm}
\label{tpert}
 Let $T>0$ and let $0<s<4$. Let $(\rhot,\vt,\ct)$ be the constant state in (\ref{constantstates}). There exists constants $\delta,\e_0\in
(0,1)$ such that, if the initial data $\rho_0,v_0,c_0$ satisfy
\begin{equation*}
  \|\rho_0^\e-\rhot\|_{H^4(\mathbb{T}^2)} +
  \e\|v_0^\e\|_{H^4(\mathbb{T}^2)}+\sqrt{\e}\|c_0^\e-\ct\|_{H^4(\mathbb{T}^2)}
  \leq \delta,
\end{equation*}
for all $\e\in(0,\e_0)$, then there exists a classical solution
$(\rho^\e,v^\e, c^\e)$ to (\ref{eq:pert1}) such that the quantity
\begin{equation*}
  \sup_{0\leq t\leq
T}\left(\|\rho^\e(t)\|_{H^s(\mathbb{T}^2)}+\e
\|v^\e(t)\|_{H^s(\mathbb{T}^2)}+\sqrt{\e}\|c^\e(t)\|_{H^s(\mathbb{T}^2)}\right)
\end{equation*}
is uniformly bounded with respect to $\e\in(0,\e_0)$ and such that
the density $\rho^\e$ satisfies
\begin{equation*}
\rho^\e (x,t)>\rhot/2>0
\end{equation*}
for all $\e\in(0,\e_0)$.
\end{thm}

\proof For any fixed $\e\in(0,\e_0)$, the sequence
$(\rho^n,v^n,c^n)$ has all space derivatives up to order $4$ in
$L^2$ and all time derivatives up to order $3$ in $L^2$.
Therefore, $(\rho^n,v^n,c^n)$ is relatively strongly compact in
$W^{1,\infty}$ and it converge (up to a subsequence) to a solution
to the original problem (\ref{eq:pert1}). Moreover, the estimate
\begin{equation*}
  \sup_{0\leq t\leq
T}\left(\|\rho^\e(t)\|_{H^s(\mathbb{T}^2)}+\e
\|v^\e(t)\|_{H^s(\mathbb{T}^2)}+\sqrt{\e}\|c^\e(t)\|_{H^s(\mathbb{T}^2)}\right)\leq
K
\end{equation*}
can be passed to the limit by weak lower semicontinuity and the
proof is complete.
\endproof

\begin{rem}
\emph{The whole procedure developed in the proof of the above
theorem can be easily generalized to the case of the third scaling
introduced in section \ref{sec:thirdscaling}.}
\end{rem}

\begin{rem}\label{remblowup}\emph{We observe here that the power like expression for the
pressure $g(\rho)=\rho^\gamma$ can be replaced by a more general
one in order to achieve the same existence result as in the above
theorem. In particular one can use $g(\rho)=\log \rho$, thus
obtaining a system which relaxes toward a Keller--Segel type
system with linear diffusion. Therefore, some of the relaxation
results contained in chapter \ref{chapconv} would include
Keller--Segel type system with linear diffusion as possible
limits. This fact is not in contradiction with the finite time
blow up phenomena occurring in the latter, because the class of
initial data for which the above theorem holds is not significant
enough in order to see the appearance of blow--up in the limit
system.}
\end{rem}

\section*{Acknowledgements}
The authors are partially supported by the Italian GNAMPA 2006
research project `Sistemi iperbolici nonlineari: modelli e analisi
qualitativa' of the `Istituto Nazionale di Alta Matematica -
Gruppo Nazionale per l'Analisi Matematica, la Probabilit\`{a} e le
loro Applicazioni'. MDF acknowledges support from the Wolfgang
Pauli Institute of Vienna and from the Wittgenstein 2000 Award of
Peter A. Markowich.

\bibliography{persistence}
\bibliographystyle{amsalpha}

\end{document}